\documentclass[11pt]{article}
\usepackage{epic}                             
\usepackage{eepic}
\usepackage{latexsym}

\begin{document}
\begin{center}\Large\bf
{
Numerical Method for Solving Obstacle Scattering Problems
by an Algorithm Based on the Modified Rayleigh Conjecture 
}
\end{center}
\begin{center}\large Weidong Chen and Alexander Ramm, \end{center}
\begin{center}\large Math Dept, Kansas State University, Manhattan, KS
66506 e-mail:chenw@math.ksu.edu and ramm@math.ksu.edu \end{center} 

{\bf Abstract.} In this paper we present a numerical algorithm for solving 
the direct scattering problems by the Modified Rayleigh Conjecture Method 
 (MRC)
introduced in [1]. Some numerical examples are given. They show that the
method is numerically efficient.



{\bf\Large I. Introduction}

 The classical Rayleigh Conjecture is discussed in [4] and [5], where
it is shown that, in general, this conjecture is incorrect: there are
obstacles (for example, sufficiently elongated ellipsoids) for
which the series, representing the scattered field outside a ball 
containing the obstacle, does not converge up to the boundary of this 
obstacle.

 The Modified Rayleigh  Conjecture (MRC) has been formulated and proved in 
[1] (see Theorem 1 below). A numerical method for solving obstacle 
scattering problems, based on MRC, was proposed in [1]. This method was 
implemented in [2] for two-dimensional
obstacle scattering problems. The numerical results in [2] were quite
 encouraging: they show that the method is efficient, economical, and is
quite competitive compared with the usual boundary integral equations
method (BIEM). A recent paper [3] contains a numerical implementation of
MRC in some three-dimensional obstacle scattering problems. Its
results reconfirm the practical efficiency of the MRC method.
   
 In this paper a numerical implementation of the Modified
Rayleigh Conjecture (MRC) method  for solving obstacle
scattering problem in three -dimensional case is
presented. Our aim
is to consider more general than in [3] three-dimensional 
obstacles: non-convex, non-starshaped, non-smooth, and to study the 
performance of the MRC in these cases. The minimization problem (5) (see 
below),
which is at the heart of the MRC method, is treated numerically in a new 
way, different from the one used in [2] and [3]. 
Our results present further numerical evidence of the 
practical efficiency of the MRC method
for solving obstacle scattering problems.

 The obstacle scattering
problems (1)-(3), we are interested in, consists of solving the equation 
$$
~~~~~~~(\bigtriangledown ^2 + k^2)u=0 ~~~~ in~~ D'=R^3\setminus D,
~~~~~~~~~~~~~~~~~~~~~~~(1) $$ where $D\subset R^3$ is a bounded domain,
satisfies the Dirichlet boundary condition $$ ~~~~~~~~~~~~~~~~~~~u|_S
=0,~~~~~~~~~~~~~~~~~~~~~~~~~~~~~~~~~~~~~~~~~~~~~~~(2) $$ where $S$ is the
boundary of $D$, which is assumed Lipschitz in this paper, and the
radiation condition at infinity: $$ u=u_0 +v=u_0 +A(\alpha ', \alpha )
\frac {e^{ikr}}{r} + o(\frac 1 r)~~~\quad r\rightarrow \infty ,~~~~~~~(3) 
$$ 
$$ r:=|x|,~~
\alpha '=x/r, ~~u_0:=e^{ik\alpha\cdot x}, $$ where $v$ is the scattered
field, $\alpha\in S^2$ is given, $S^2$ is the unit sphere in $R^3$,
$k=const>0$ is fixed, $k$ is the wave number. The coefficient $A(\alpha
',\alpha)$ is called the scattering amplitude.

Denote
$$
~~~~~~~~~~~~~A_l (\alpha):=\int_{S^2} A(\alpha ', \alpha) 
\overline {Y_l(\alpha ')}d\alpha ',\qquad~~~~~~~~~~~~~~~~~~~~~~~~(4)
$$
where $Y_l(\alpha )$ are the orthonormal spherical harmonics:
$$
Y_l=Y_{lm},~~ -l\leq m \leq l,~~l=0,1,2,...
$$
$$
Y_{lm}(\theta,\phi)=\frac 1 {\sqrt {2\pi}} e^{im\phi} \Theta _{lm}(cos\theta),
$$
$$
\Theta _{lm}(x) =\sqrt {\frac {2l+1} 2  \frac {(l-m)!} {(l+m)!}} P^m_l(x),
$$
$P^m_l(x)$ are the associated Legendre functions of the first kind,
$$
P^m_l(x):=(1-x^2)^{m/2} \frac {d^mP_l(x)} {dx^m},~~m\geq 0, 
$$
and
$$
P_l(x):=\frac {(-1)^l} {2^l l!} \frac {d^l} {dx^l} (1-x^2)^l.
$$
For $m<0$
$$
\Theta _{lm}(x) =(-1)^m \Theta _{l,-m}(x).
$$
Let $h_l(r)$ be the spherical Hankel functions of the first kind, normalized so that $h_l(kr) \sim e^{ikr}/r$ as $r \rightarrow +\infty$. Let $B_R:=\{x:|x|\leq R\}\supset D$, and the origin is inside $D$.

Then in the region $r>R$, the solution to the acoustic wave problem (1)-(3) is of the form:
$$
u(x, \alpha)=e^{ik\alpha\cdot x} + \sum_{l=0}^{\infty} A_l(\alpha) \psi _l(x),~|x|>R,
$$
$$
\psi _l:=Y_l(\alpha ')h_l(kr),~~~~~ r>R,~~~\alpha '=x/r,
$$
where 
$$
\sum_{l=0}^{\infty}:=\sum_{l=0}^{\infty}\sum_{m=-l}^{l}.
$$
Fix $\epsilon >0$, an arbitrary small number. The following Lemmas and Theorem1 are proved in [1]. 

 {\bf Lemma 1.} {\em There exist $L=L(\epsilon )$ and numbers $c_l=c_l(\epsilon )$ such that}
 $$
~~~~~~~~~~~~~ ||u_0(s)+\sum_{l=0}^{L} c_l(\epsilon ) \psi _l(s)||_{L^2(S)}< \epsilon.~~~~~~~~~~~~~~~~~~~~~~~(5)
 $$
 
 {\bf Lemma 2.}{\em If (5) holds, then $||v_{\epsilon}(x)-v(x)||=O(\epsilon),~\forall x\in D',~~~ \epsilon\rightarrow 0.$
 where
 $$
 ~~~~~~~~~~~~~~~v_\epsilon (x):=\sum_{l=0}^{L} c_l(\epsilon ) \psi _l(x),~~~~~x\in D',~~~~~~~~~~~~~~~~~~~(6)
 $$
 and
 $$
 ||.||:=||.||_{H^m_{loc}(D')} +||.||_{L^2(D'; (1+|x|)^{-\gamma)}},~\gamma >0, m>0,~~~~~~~~~~~(7)
 $$
$m$ is arbitrary, and $H^m$ is the Sobolev space}.

{\bf Lemma 3.} $c_l(\epsilon) \rightarrow A_l(\alpha), \forall l, \epsilon \rightarrow 0.$
 
 {\bf Theorem 1 (Modified Rayleigh Conjecture).}
{\em Let $D\in R^3$ be a bounded obstacle with Lipschitz boundary
S. For any $\epsilon>0$ there
exists $L=L(\epsilon)$ and $c_l(\epsilon)=c_{lm}(\epsilon)$, $0\leq l\leq L$, $-l\leq m\leq l$, such that inequality (5) holds. If (5) holds
then function (6) satisfies the estimate
$||v(x)-v_\epsilon(x)||=O(\epsilon)$, where the norm is defined in (7).
Thus, $v_\epsilon(x)$ is an approximation of the scttered field
everywhere in $D'$.}

 In order to obtain an accurate solution, usually one has to take $L$ large. But as $L$ grows the condition number of the matrix $(\psi _l,~\psi _{l'})_{L^2(S)}$ is increasing very fast. So we choose some interior points $x_j\in D, ~j=1,2,...,J$, and use the following version of Theorem 1([2]):
 
 {\bf Theorem 2.} {\em Suppose $x_j\in D, ~j=1,2,...,J$, then $\forall\epsilon >0,~\exists L=L(\epsilon)$ and $c_{lj}(\epsilon), l=0,...,L,~j=0,...,J(\epsilon)$, such that
 
(i) $$
 ||u_0(s)+\sum_{j=0}^{J}\sum_{l=0}^{L} c_{lj}(\epsilon ) \psi _l(s-x_j)||_{L^2(S)}< \epsilon.~~~~~~~~~~~~~~~(5')
 $$
 
 (ii)$$||v_{\epsilon}(x)-v(x)||=O(\epsilon),$$
 where
 $$
 v_{\epsilon}(x)=\sum_{j=0}^{J}\sum_{l=0}^{L} c_{lj}(\epsilon ) \psi _l(s-x_j)
 $$
 and the $||.||$ is defined in Lemma 2.}
 
 {\bf Remark.} Theorem 1 is the basis for MRC algorithm for
computation of the field scattered by an obstacle: one takes
an $\epsilon>0$ and an integer $L>0$,
minimizes the left-hand side of (5) with respect to $c_l$,
and if the minimum is $\leq \epsilon$ then the function (6) is
the
approximate solution of the obstacle scattering problem
with the accuracy $O(\epsilon)$.
If the above minimum is greater than $\epsilon$, then one
increases $L$ until the minimum is less than $\epsilon$. This is
possible by Lemma 1.
In computational practice, one may increase also the number
$J$
of points $x_j$ inside $D$, as explained in Theorem 2. The
increase of $J$ allows one to reach the desired value of the
above minimum keeping $L$ relatively small. This gives
computational advantage in many cases.

 In section 2, an algorithm is presented for solving the problem (1)-(3).
 This algorithm is based on the MRC. Compared with the previous work in
the case of two- and three-dimensional MRC([2],[3]), we consider more
general surfaces, in particular non-starshaped and piecewise-smooth
boundaries. The numerical results are given in section 3. A discussion of
the numerical results is given in section 4.
  
{\bf\Large II. The MRC algorithm for Solving Obstacle Scattering Problems}

{\bf 1. Smooth starshaped boundary:}

Assume the surface $S$ is given by the equation
$$
r =r (\theta , \varphi),~~~~~~~~~~~ 0\leq \varphi\leq 2\pi,~0\leq \theta\leq \pi.
$$
Define
$$
~~~~~~~~~~F(c_0, c_1,...,c_L):=||u_0+\sum_{l=0}^L c_l \psi_l ||_{L^2(S)}^2.~~~~~~~~~~~~~~~~~~~~~~~(5'')
$$
Let
$$
h_1=2\pi/n_1,~~~~~~h_2=\pi/n_2
$$
$$
0=\varphi_0<\varphi_1<....<\varphi_{n_1}=2\pi,~ \varphi_{i_1}=i_1h_1,~~i_1=1,...,n_1,
$$
$$
0=\theta_0<\theta_1<....<\theta_{n_1}=\pi,~ \theta_{i_2}=i_2h_2,~~i_2=1,...,n_2,
$$
where $n_1$ and $n_2$ are the number of steps.
By Simpson's formula([8]), we obtain an approximation of $F(c_0, c_1,...,c_L)$:
$$
F(c_0, c_1,...,c_L)=\sum_{i_1=0}^{n_1} \sum_{i_2=0}^{n_2} a_{i_1 i_2}  {\huge |}u_{0i_1i_2}+\sum_{l=0}^L c_l \psi_{l i_1 i_2}{\huge |}^2 w_{i_1 i_2} h_1 h_2~~~~~~~(5''')
$$
where
\vphantom{construct}
\[ a_{i_1, i_2}=
   \left\{ \begin{array}{l}
~4,~~~~~i_1~ and~ i_2~ even 
\\~8,~~~~~i_1~ -~i_2~ odd 
\\16,~~~~~i_1~ and~ i_2~ odd
\end{array} \right. \]
and
$$
\psi_{li_1i_2}=Y_l(\theta _{i_1}, \varphi _{i_2})h_l(kr(\theta _{i_1}, \varphi _{i_2})),~~~
w_{i_1 i_2}=w(\theta_{i_1} ,\varphi_{i_2})
$$
where
$$
~~~~~~~~~~w(\theta ,\varphi )=(r^2 r^2_\varphi + r^2 r^2_\theta sin^2 \theta +r^4 sin^2 \theta )^{1/2}.~~~~~~~~~~~~~~~~~~~~~(8)
$$
We can find  $c^*=(c_0^*, c_1^*,...,c_L^*)$ such that
$$
~~~~~~~~~~~~~~~~~~~F(c^*)=min F(c_0, c_1,...,c_L).~~~~~~~~~~~~~~~~~~~~~~~~~~~~~~~(9)
$$
We first write
$$
~~~~~~~~~~~~~~~~~~~~~~~~~~F(c)=||Ac-B||^2,~~~~~~~~~~~~~~~~~~~~~~~~~~~~~~~~~~~(10)
$$
where
$$
A=(A_{l,i})_{M\times L_1},~~A_{l,i}=\psi_{l i_1 i_2} (a_{i_1 i_2} w_{i_1 i_2} h_1 h_2)^{\frac 1 2},~~i=i_1i_2,
$$
$$
B=(B_i)_{M\times 1},~~B_i=u_{0i_1i_2} (a_{i_1 i_2} w_{i_1 i_2} h_1 h_2)^{\frac 1 2},
$$
in which $M=n_1n_2$, $L_1=(L+1)(2L+1)$ since $c_l=c_{lm}, 0\leq l\leq L, -l\leq m\leq l$.

Then Householder reflections are used to compute an orthogonal-triangular factorization: $A*P = Q*R$
where P is a permutation([8], p.171), Q is an orthogonal matrix, and R is an upper triangular matrix. Let $r=rank(A)$.
This algorithm requires $4ML_1r-2r^2(M+L_1)+4r^3/3$ flops([9], pp.248-250).
The least squares solution $c$ is computed by the formula $c = P*(R^{-1}*(Q'*(A^TB)))$. This minimization procedure is based on the matlab code([10]).

In [2] and [3] singular value decomposition was used for minimization of (5''). Here we use 
the matlab minimization code which is based on a factorization of the matrix A. This has 
the following advantages from the point of view of numerical analysis. We can choose an integer $r_1$:
$$
0<r_1\leq r
$$
such that the first $r_1$ rows and columns of $R$ form a well-conditioned matrix when A is not of full rank, or the
rank of A is in doubt([10]). See Golub and Van Loan [9] for a further discussion of numerical rank determination. 

If we choose $x_j\in D,~ j=1,...,J$, we obtain
$$
F_J(c)=F_J(c_{01},..., c_{0J}, c_{11},...,c_{1J},..., c_{L1},...,c_{LJ})
$$
$$
=\sum_{i_1=0}^{n_1} \sum_{i_2=0}^{n_2} \sum_{j=1}^{J} a_{i_1 i_2}|u_{0i_1i_2}+\sum_{l=0}^L c_{lj} \psi_{l i_1 i_2}|^2 w_{i_1 i_2} h_1 h_2.
$$
The algorithm for finding the minimum of $F_J(c)$ will be same.

{\bf 2. Piecewise-smooth boundary:}

Suppose 
$$
S=\bigcup _{n=1} ^N S_n.
$$
Then
$$
F(c_0, c_1,...,c_L)=\sum_{n=1}^N||u_0+\sum_{l=0}^L c_l \psi_l ||_{L^2(S_n)}^2
$$
$$
\forall (x,y,z) \in S_n,~~ r ^2 =x^2 +y^2 +z^2, ~~\cos\theta =z/r,~~\tan \varphi =y/x.~~~~~~(11)
$$

{\bf 3. Non-starshaped case:}

Suppose $S$ is a finite union of the surfaces, each of which is starshaped with respect to a point $\vec{r^0 _n}$,
$$
S=\bigcup _{n=1} ^N S_n.
$$
and the the surfaces $S_n$ are given by the equations in local spherical coordinates:
$$
S_n:~~~~ \vec{r}-\vec{r^0 _n}=(r_n(\theta _n , \varphi _n) \cos\varphi _n\sin\theta _n,~ r_n(\theta _n, \varphi _n) \sin\varphi _n\sin\theta _n, ~r _n(\theta _n, \varphi _n) \cos\theta _n),
$$
$$
n=1,...N,
$$
where $\vec{r^0 _n}$ are constant vectors.

Then
$$
F(c_0, c_1,...,c_L)=\sum_{n=1}^N||u_0+\sum_{l=0}^L c_l \psi_l ||_{L^2(S_n)}^2.
$$
The weight functions $w_n(\theta, \varphi)$ are the same as in (8) since $\vec{r^0 _n}$ are constant vectors.

{\bf\Large III. Numerical Results}

     In this section, we give four examples to show the convergence
rate of the algorithm and how the error depends on the shape of $S$.

{\bf Example 1.} The boundary S is the sphere of radius 1 centered at the origin.

In this example, the exact coefficients are:
$$
c_{lm}=- \frac {4\pi i^l j_l(k)} {h_l(k)}~ \overline{Y_{lm}(\alpha)}
$$
Let $k=1,~~\alpha=(1,0,0)$.
We choose $n_1=20, ~n_2=10$.

---------------------------------------------------------------------------------------------
  
  L~~~~~~~~0~~~~~~~~1.0000~~2.0000~~3.0000~~4.0000~~5.0000~~6.0000~~7.0000

--------------------------------------------------------------------------------------------

$F(c^*)$   ~6.3219~~1.6547~~0.2785~~0.0368~~0.0034~~0.0003~~0.0000~~0.0000

--------------------------------------------------------------------------------------------

err($c$)   ~0.0303~~0.0172~~0.0020~~0.0004~~0.0000~~0.0000~~0.0000~~0.0000

--------------------------------------------------------------------------------------------

where $$err(c)=(\sum_{l=0}^L |c^*_l-c_l|^2)^{\frac 1 2}.$$
When $n_1=40, ~n_2=20,$

---------------------------------------------------------------------------------------------
  
  L~~~~~~~~0~~~~~~~~1.0000~~2.0000~~3.0000~~4.0000~~5.0000~~6.0000~~7.0000

--------------------------------------------------------------------------------------------

$F(c^*)$   ~6.3544 ~~1.6562~~0.2820~~0.0358~~0.0036~~0.0003~~0.0000~~0.0000

--------------------------------------------------------------------------------------------

err($c$)   ~0.0147~~0.0076~~0.0011~~0.0001~~0.0000~~0.0000~~0.0000~~0.0000

--------------------------------------------------------------------------------------------

Next, we fix $n_1=20, ~n_2=10$ and test the  results for different k and $\alpha$.

When $k=2,~~\alpha=(1,0,0)$,

---------------------------------------------------------------------------------------------
  
  L~~~~~~~~0~~~~~~~~1.0000~~2.0000~~3.0000~~4.0000~~5.0000~~6.0000~~7.0000

--------------------------------------------------------------------------------------------

$F(c^*)$   ~10.4506~~5.5783~~1.9291~~0.5217~~0.0970~~0.0156~~0.0020~~0.0003

--------------------------------------------------------------------------------------------

err($c$)   ~0.0404~~0.0205~~0.0048~~0.0020~~0.0005~~0.0000~~0.0000~~0.0000

--------------------------------------------------------------------------------------------

When $k=1,~~\alpha=(0,1,0)$,

---------------------------------------------------------------------------------------------
  
  L~~~~~~~~0~~~~~~~~1.0000~~2.0000~~3.0000~~4.0000~~5.0000~~6.0000~~7.0000

--------------------------------------------------------------------------------------------

$F(c^*)$   ~6.3801~~1.6628~~0.2821~~0.0371~~0.0044~~0.0003~~0.0000~~0.0000

--------------------------------------------------------------------------------------------

err($c$)   ~0.0014~~0.0106~~0.0005~~0.0004~~0.0000~~0.0000~~0.0000~~0.0000

--------------------------------------------------------------------------------------------

When $k=1,~~\alpha=(0,0,1)$,

---------------------------------------------------------------------------------------------
  
  L~~~~~~~~0~~~~~~~~1.0000~~2.0000~~3.0000~~4.0000~~5.0000~~6.0000~~7.0000

--------------------------------------------------------------------------------------------

$F(c^*)$   ~6.4156~~1.6909~~0.2955~~0.0418~~0.0025~~0.0002~~0.0000~~0.0000

--------------------------------------------------------------------------------------------

err($c$)   ~0.0093~~0.0109~~0.0049~~0.0007~~0.0001~~0.0000~~0.0000~~0.0000

--------------------------------------------------------------------------------------------

When $k=1,~~\alpha=(1/\sqrt{2},1/\sqrt{2},0)$,

---------------------------------------------------------------------------------------------
  
  L~~~~~~~~0~~~~~~~~1.0000~~2.0000~~3.0000~~4.0000~~5.0000~~6.0000~~7.0000

--------------------------------------------------------------------------------------------

$F(c^*)$   ~6.3500~~1.6711~~0.2810~~0.0371~~0.0040~~0.0003~~0.0000~~0.0000

--------------------------------------------------------------------------------------------

err($c$)   ~0.0218~~0.0057~~0.0019~~0.0004~~0.0001~~0.0000~~0.0000~~0.0000

--------------------------------------------------------------------------------------------

When $k=1,~~\alpha=(1/\sqrt{3},1/\sqrt{3},1/\sqrt{3})$,

---------------------------------------------------------------------------------------------
  
  L~~~~~~~~0~~~~~~~~1.0000~~2.0000~~3.0000~~4.0000~~5.0000~~6.0000~~7.0000

--------------------------------------------------------------------------------------------

$F(c^*)$   ~6.3739~~1.6542~~0.2850~~0.0368~~0.0040~~0.0003~~0.0000~~0.0000

--------------------------------------------------------------------------------------------

err($c$)   ~0.0170~~0.0054~~0.0021~~0.0003~~0.0001~~0.0000~~0.0000~~0.0000

--------------------------------------------------------------------------------------------

{\bf Example 2.} The boundary S is the surface of the cube $[-1,1]^3$.
Here
$$
S=\bigcup _{n=1} ^6 S_n.
$$
and
$$
F(c_0, c_1,...,c_L)=\sum_{n=1}^6||u_0+\sum_{l=0}^L c_l \psi_l ||_{L^2(S_n)}^2
$$
$$
=\sum_{n=1}^6\sum_{i_1=0}^{n_1} \sum_{i_2=0}^{n_2} a_{i_1 i_2}  {\huge |}
u_{0i_1i_2}+\sum_{l=0}^L c_l \psi_{l i_1 i_2}{\huge |}^2 \Delta_1 \Delta_2
$$
where
$$
\Delta_1=2/n_1,~~~~~~\Delta_2=2/n_2.
$$
The origin is chosen at the cenetr of symmetry of the cube. The surface area 
element is calculated in the Cartesian coordinates, so the weight $w=1$.

Let $S_1$ be the surface
$$
z=1, ~-1\leq x\leq 1, ~-1\leq y\leq 1
$$ 
and
$$
x_{i_1}=-1+i_1 \Delta_1,~0\leq i_1\leq n_1
$$
$$
y_{i_2}=-1+i_2 \Delta_2,~0\leq i_2\leq n_2
$$
Then
$$
\psi_{li_1i_2}=Y_l(\theta _{i_1}, 
\varphi _{i_2})h_l(kr(\theta _{i_1}, \varphi _{i_2})),~~~
$$
and $\theta _{i_1}$ and $\varphi _{i_2}$ can be computed by formula (11). 
For other surfaces $S_j$ the algorithm is similar. 

The values of $\min F(c)=F(c*)$ and the values $\min F_J(c)=F_J(c*)$ with $x_j$:
$$
\{x_j:j=0,...,6\}=\{(0,0,0), (0.2,0,0), (-0.2,0,0), 
$$
$$
(0,0.2,0), (0,-0.2,0), (0,0,0.2),(0,0,-0.2)\}
$$
are given below.

We choose $n_1=10, ~n_2=10$

-----------------------------------------------------------------------------------------------
  
  L~~~~~~~~~~~~~~0~~~~~~~~~~1.0000~~2.0000~~3.0000~~4.0000~~5.0000~~6.0000~~7.0000~~8.0000

---------------------------------------------------------------------------------------------

$F(c^*)$~~~~10.6301~~3.6277~~2.6760~~2.2309~~1.9832~~1.5737~~1.5034~~1.2948~~1.1753

---------------------------------------------------------------------------------------------

$F_J(c^*)$~~~2.6297~~1.0970~~0.5487~~0.1572~~0.0667~~0.0320~~0.0168~~0.0078~~0.0035

---------------------------------------------------------------------------------------------

When $n_1=20, ~n_2=20,$

-----------------------------------------------------------------------------------------------
  
  L~~~~~~~~~~~~~~0~~~~~~~~~~1.0000~~2.0000~~3.0000~~4.0000~~5.0000~~6.0000~~7.0000~~8.0000

---------------------------------------------------------------------------------------------

$F(c^*)$~~~~10.7923~~3.7144~~2.7778~~2.3393~~2.0873~~1.6671~~1.5938~~1.4277~~1.3368

---------------------------------------------------------------------------------------------

$F_J(c^*)$~~~2.7248~~ 1.1433~~0.5757~~0.1686~~0.0694~~0.0652~~0.0236~~0.0143~~0.0090

---------------------------------------------------------------------------------------------

{\bf Example 3.} The boundary S is the surface of the ellipsoid
 $x^2+y^2+z^2/b^2=1,$ the values of $\min F(c)=F(c*), ~ b=2,3,4,5$ 
with $n_1=20,~n_2=10$ are:

---------------------------------------------------------------------------------------------

L~~~~~~~~0~~~~~~~~~1.0000~~2.0000~~3.0000~~4.0000~~5.0000~~6.0000~~7.0000

---------------------------------------------------------------------------------------------

b=2~~8.8836~~5.4955~~3.0421~~2.8434~~1.3622~~ 1.2093~~0.8753~~0.8132

---------------------------------------------------------------------------------------------

b=3~~14.1617~~12.0477~~7.2296~~7.0999~~3.8077~~3.6829~~3.1324~~3.0496

---------------------------------------------------------------------------------------------

b=4~~19.5326~~17.9346~~9.9927~~9.8720~~5.3333~~5.2008~~4.6793~~4.5738

----------------------------------------------------------------------------------------------

b=5~~22.9765~~21.5653~~11.4850~~11.3587~~6.1637~~6.0096~~5.5202~~5.3933

----------------------------------------------------------------------------------------------

The values of $\min F_J(c)=F_J(c*), ~ b=2,3,4,5$ with $x_j$:
$$
\{x_j:j=0,...,6\}=\{(0,0,0), (0.5,0,0), (-0.5,0,0),
$$
$$
(0,0.5,0), (0,-0.5,0), (0,0,0.5),(0,0,-0.5)\}
$$
are:

---------------------------------------------------------------------------------------------

L~~~~~~~~0~~~~~~~~~1.0000~~2.0000~~3.0000~~4.0000~~5.0000~~6.0000~~7.0000

---------------------------------------------------------------------------------------------

b=2~~ 2.4856~~0.7090~~0.2530~~0.0062~~0.0000~~0.0000~~0.0000~~0.0000

---------------------------------------------------------------------------------------------

b=3~~4.6639~~1.3619~~0.6618~~0.0074~~0.0000~~0.0000~~0.0000~~0.0000

---------------------------------------------------------------------------------------------

b=4~~5.5183~~1.8624~~0.7844~~0.0060~~ 0.0000~~ 0.0000~~ 0.0000~~ 0.0000

----------------------------------------------------------------------------------------------

b=5~~11.0579~~8.7027~~6.4831~~0.8357~~0.0017~~0.0000~~0.0000~~0.0000

----------------------------------------------------------------------------------------------

{\bf Example 4.} The obstacle is a dumbbell. Its boundary S is not 
smooth, non-starshaped and not convex:
$$
S=S_1\bigcup S_2\bigcup S_3
$$
$$
S_1:~\vec{r} -(0,0,1)=(1.5\cos\varphi\sin\theta, 1.5\sin\varphi\sin\theta, 1.5\cos\theta )
$$
$$
S_2:~\vec{r} -(0,0,-1)=(1.5\cos\varphi\sin\theta, 1.5\sin\varphi\sin\theta, 1.5\cos\theta )
$$
$$
S_3:r \sin\theta =1
$$
$$
\{x_j:j=0,...,10\}=\{(0, 0, 0), (0, 0, 0.1), (0, 0, -0.1), (0, 0, 0.2),(0, 0, -0.2), 
$$
$$
(0, 0, 0.3), (0, 0, -0.3), (0,0,0.4),(0,0,-0.4), (0,0,0.5),(0,0,-0.5)\};
$$
We choose $n_1=20$, $n_2=10$ for every $S_i(i=1,2,3)$.

-------------------------------------------------------------------------------------------

L~~~~~~~~~~~~0~~~~~~~~~~~1.0000~~~2.0000~~~3.0000~~~4.0000~~~5.0000~~~6.0000~~~7.0000

-------------------------------------------------------------------------------------------

$F(c^*)$~~~~~~25.8840~~ 20.8059~~16.4968~~15.6622~~12.9241~~12.1915~~11.0187~~9.5263

-------------------------------------------------------------------------------------------

$F_J(c^*)$~~~~20.3118~~~  8.0238~~~~5.1062~~~2.5908~~~0.8304~~~0.4067~~~0.0453~~~~0.0084

-------------------------------------------------------------------------------------------

{\bf\Large IV. Conclusion}

From the numerical results one can see that the accuracy of the numerical
solution depends on the smoothness and elongation of the object.

In Example 1 the surface $S$ is a unit sphere and the numerical solution 
is very accurate. In Example 3 the results for different elongated 
ellipsoids
show that if the elongation (eccentricity) grows, then the accuracy 
decreases. In Example 2 the surface
is not smooth and  the result is less accurate than in Example 3. In 
Example 4 the surface in nonconvex and not smooth, but the accuracy is 
of the same order as in Example 2.

When $b$ is large or $S$ is not smooth, the numerical results in Example 2
and Example 3 show that if one adds more points $x_j$ then the accuracy
of the solution increases.

In Example 1 and Example 2, as one increased $n_1$ and $n_2$, the minimum 
$F(c*)$ has also increased because the condition number of the matrix $A$ 
in (10) grew as $n_1$ and $n_2$ increased.

Using the results of Example 1 one can check the accuracy
in finding $c_l$ by the value of the minimum $$ F(c*)\leq \epsilon. $$ 
\newpage
\begin{center}\Large\bf
                              References
\end{center}

[1] Ramm A. G. [2002], Modified Rayleigh Conjecture and Applications, J.
Phys. A: Math. Gen. 35, L357-L361.

[2] Gutman S. and Ramm A. G. [2002], Numerical Implementation of the MRC
Method for Obstacle Scattering Problems, J. Phys. A: Math. Gen. 35,
L8065-L8074.

[3] Gutman S. and Ramm A. G., Modified Rayleigh Conjecture Method for
Multidimensional Obstacle Scattering Problems(submitted).

[4] Barantsev, R., Concerning the Rayleigh hypothesis in the problem of
scattering from finite bodies of arbitrary shapes, Vestnik Lenigrad.
Univ., Math., Mech., Astron., 7, (1971), 52-62.

[5] Millar, R., The Rayleigh hypothesis and a related least-squares
solution to scattering problems for periodic surfaces and other
scatterers, Radio, Sci., 8, (1973), 785-796.

[6] Ramm, A. G., Scattering by obstacles, D. Reidel, 1986.

[7] Triebel H., Theory of Function Spaces, vol. 78 of Monographs in
Mathematics. Birkhauser Verlag, Basel, 1983.

[8] Kincaid D. and Cheney W., Numerical Analysis: Mathematics of
Scientific Computing, Brooks/Cole, 2002.

[9] Golub G. H. and Van Loan C. F., Matrix Computations, The John Hopkins University Press: Baltimore and London, 1996.

[10] Anderson, E., Z. Bai, C. Bischof, S. Blackford, J. Demmel, J. Dongarra, J. Du Croz, 
    A. Greenbaum, S. Hammarling, A. McKenney, and D. Sorensen, LAPACK User's Guide 
    (http://www.netlib.org/lapack/lug/lapack\_lug.html), Third Edition, SIAM, Philadelphia, 1999.  
\end{document}